\documentstyle[amstex,amssymb]{article}

\def\m@th{\mathsurround=0pt }
\def\cases#1{\left\{\,
\vcenter{\normalbaselines
\m@th\ialign{$##\hfil$&\quad##\hfil\crcr#1\crcr}}\right.}
\newcommand{\th}{{\theta}}
\def\de{{\delta}}

\def\al{{\alpha}}

\def\be{{\beta}}

\def\de{{\delta}}

\newcommand{\bbb}{{\cal B}}

\newcommand{\Ee}{{\bold E}}
\newcommand{\Ff}{{\bold F}}
\newcommand{\Kk}{{\bold K}}

\newcommand{\PP}{{\Bbb P^1}}

\newcommand{\F}{{\Bbb F}}
\newcommand{\K}{{\Bbb K}}
\newcommand{\ab}{\hspace{10mm}}

\newcommand{\G}{{\bold G}}

\newcommand{\HH}{{\bold H}}

\newcommand{\mm}{{\cal M}}

\newcommand{\UU}{{\bold U}_{\bold q}}
\newcommand{\Uu}{{\bold U}}

\newcommand{\Oc}{{\cal O}}

\newcommand{\g}{{\goth g}}

\newcommand{\C}{{\Bbb C}}

\newcommand{\Z}{{\Bbb Z}}

\def\L{{\Bbb L}}

\def\End{{\mbox{End}}}

\newcommand{\gt}{{\g_{_{tor}}}}
\def\End{{\mbox{End}\,}}

\newcommand{\zw}{{{\textstyle{\small \frac{z}{w}}}}}

\begin{document}

\setlength{\parindent}{0pt}
\setlength{\parskip}{3pt plus 5pt minus 0pt}

\centerline{\large{\bf LANGLANDS RECIPROCITY FOR ALGEBRAIC SURFACES}}
\vskip 7mm
\centerline{\sc VICTOR GINZBURG , MIKHAIL KAPRANOV , ERIC VASSEROT}

\vskip 10mm

\ab Twenty-five years ago R. Langlands proposed [L]
a ``fantastic generalization" of Artin-Hasse
reciprocity law in the classical class field theory. He
conjectured  the existence  of a correspondence
between automorphic  irreducible
infinite-dimensional representations of a reductive group $G$
over a global number field on the one hand,
and (roughly speaking) finite dimensionsional representations
of the Galois group of the field, on the other hand.
Following an earlier idea of A.Weil, Langlands' conjecture
was reinterpreted by V.Drinfeld [Dr 3] and G. Laumon [La]
in purely geometric terms. In the complex geometry setup the number field
gets replaced by a Riemann surface  $X$, the Galois group
of the field gets replaced by the fundamental group, $\pi_1(X)$, and
an automorphic representation gets replaced by a perverse sheaf, cf. [BBD],
on the moduli space of algebraic principal
$G$-bundles on $X$.

\ab The purpose of this note is to formulate some results and
conjectures
related to
Langlands reciprocity for algebraic surfaces, as opposed to  algebraic
curves. That would amount, in number theory, to
a non-abelian higher dimensional
generalization of the class field theory, the subject
that remains a total mystery at the moment. A key step of our
geometric approach is a construction of Hecke operators
for  vector
bundles on an algebraic surface. The main point
of this paper is that in certain cases
 the corresponding algebra of Hecke operators turns out
to be a homomorphic image of the {\it quantum toroidal algebra}.
The latter is a quantization, in the spirit of Drinfeld
[Dr 1] and Jimbo [Ji], of the universal enveloping algebra of the
universal
central extension of a ``double-loop" Lie algebra.
This yields, in particular, a new geometric
construction of affine quantum groups of types $A^{(1)}, D^{(1)}$ and
$E^{(1)}$ in terms of Hecke operators
for an elliptic surface.

\ab Our approach was motivated in part by a relation between
instantons on ALE-spaces and affine Lie algebras discovered by
Nakajima [Na 1] (which was in its turn motivated by Lusztig's
construction [Lu] of quantum groups in terms of quivers and the
{\it Lagrangian construction} of [Gi]). Nakajima's results were also
used by Vafa and Witten [VW] in verifying a special case of
their $S$-duality conjecture. We believe, that there is in fact
a very close connection between
the Langlands' reciprocity for an algebraic surface and
the $S$-duality conjecture for the underlying real 4-manifold.
Suffices it to say
that the interplay between a reductive group $G$ and the Langlands
dual group $G^\vee$ is quite essential in
both cases.

\ab Proofs of the results announced in this paper will appear elsewhere.

{\sc Acknowledgements.}
The first author is pleased to thank D. Kazhdan for a discussion a few
years ago that has lead eventually to the $K$-theoretic construction
of $\S 5$.
We are especially grateful to H. Nakajima for providing us with a
revised
unpublished version of [Na2]. Theorems 2.4 and 4.2 below are based on
the same idea as the results of [Na2] for ALE-spaces.
The second author was partially supported by an NSF grant and by
Alfred P. Sloan Research Fellowship.

{\newcommand{\Mgr}{{\Bbb C}^* \times G}
\newcommand{\BAF}{{\cal B}}
\newcommand{\RRA}{{{\bold R}(A)}}
\newcommand{\cdz}{{\Bbb C}^d((z))}
\newcommand{\FP}{{\Bbb F}}
\newcommand{\kp}{{\Bbb K}}
\newcommand{\gk}{G^\vee({\Bbb K})}
\newcommand{\lgg}{G^\vee({\Bbb K})}
\newcommand{\BAFd}{{\cal B}_{\bf d}}
\newcommand{\fd}{{\cal F}}
\newcommand{\fbaf}{{\Bbb C}_{G^\vee({\Bbb K})}[{\cal B} \times {\cal B}]}

\vskip 3mm
\pagebreak[3]
\ab {\bf 1. Hecke operators and Satake isomorphism.}

\ab We recall here the well-known,
see e.g. [Dr 3], interpretation of Hecke operators
in terms of moduli spaces.

\ab Let $C$ be a complete smooth curve over a finite field $\F_p$
and $x\in C$ an $\F_p$-rational point.
Let $\mm$ be the discrete set of
isomorphism classes of
 $\F_p$-rational vector bundles on $C$ of rank $d$. Given a vector
bundle $V\in \mm$ and $W\subset V_x$, a  linear subspace in the
fiber
of $V$ at $x$, let ${\cal A}_x(V,W)$
denote the subsheaf of
$V$ consisting of sections whose value at $x$ is contained in $W$.
This sheaf is again locally free, hence, can be regarded as a new
vector bundle on $C$, called the elementary transformation
of $V$.
Let $\C[\mm]$ be the vector space of $\C$-valued functions on $\mm$
 with finite support. For any $k=1,\ldots,d,$ the following sum is
finite,
hence defines a
Hecke operator $T_{x,k}: \C[\mm]\to \C[\mm]$ by
$$T_{x, k} (f)(V)  = \sum_{\{W\subset V_x\;|\; \dim W = d-k\}}
 f({\cal A}_x(V,W)).$$
Let $\HH_x$ denote the $\C$-subalgebra of linear endomorphisms of the
vector space $\C[\mm]$ generated by the operators
$T_{x, k}\,,\,k=1,\ldots,d.$

\ab
Let
$\C[G]^{^G}$ denote the algebra of the regular class functions on
the complex algebraic group $G=GL_d(\C)$.
The classical  theorem of Satake, combined with the fact that
the natural $\HH_x$-action makes $\C[\mm]$ a faithful $\HH_x$-module,
yields the following version of the Satake isomorphism.
\smallskip

{\sc Theorem 1.1.} {\it There is a  natural algebra
isomorphism}:
\[\C[G]^{^G}\simeq \HH_x\;.\qquad\square\]
\vskip 3mm

\ab {\bf 2. Hecke operators for algebraic surfaces.}

\ab We construct here certain analogues of Hecke operators for
 algebraic surfaces. These operators are associated to curves rather
than to
points on the surface.

\ab Let $S$ be a  smooth projective surface over a finite field
$\F_p$, and $C\subset S$ a curve, i.e., possibly non-reduced subscheme,
locally given by one equation defined over $\F_p$.
Write $\mm$ for the discrete set of the isomorphism classes
of $d$-dimensional algebraic vector bundles on $S$, i.e., of locally free
$\cal O_S$-sheaves in the Zariski topology. Given
$W\subset V|_{_C}$, a
subbundle (= subsheaf which is locally a direct summand) $W$
in the restriction of $V$ to $C$,
 we introduce the subsheaf
${\cal A}_C(V,W)$ in $V$ to be
formed by the sections of $s\in V$ such that $s|_{_C}\in W$. This sheaf is
 locally free and, hence, can be regarded as a new vector bundle on $S$,
see [Ma].

\ab Let ${\cal O}_S(-C)$ be
the invertible sheaf on $S$ formed by functions vanishing at $C$, let
${\cal O}_S(C)$ be its dual, and
${\cal O}_C(-C)$ the (scheme theoretic) restriction of ${\cal O}_S(-C)$
to $C$.
Let $\C[\mm]$
denote the vector space of complex functions on $\mm$ with finite support.
To any vector bundle $L$ on $C$ we associate
two Hecke operators $T_{_{L}},T^*_{_{L}} :$
$\C[\mm] \to \C[\mm]$ given by the following formulas
$$(T_{_{L}}f)(V) = \sum_{\{W \subset V|_{_C}\; {\small {\mbox{such that}}}
\; (V|_{_C})/W\simeq L\}}
 f({\cal A}_C(V,W))$$
$$(T^*_{_L}f)(V) =
 \sum_{\{W\subset V|_{_C}\;  {\small {\mbox{such that}}}\; W \simeq
L\otimes
{\cal O}_{_C}(-C)\}}
 f({\cal A}_C(V,W)\otimes {\cal O}_S(C)),$$

\ab Assume now that $C\simeq \PP$ is a projective
line.
Let $\mm_C\subset \mm$ be the set of those vector bundles $V$
whose restrictions to $C$ have the form
$$ V|_{_C}\simeq \cal O_{\!_\PP}(1)^{^{\oplus a(V)}}
\oplus\cal O_{\!_\PP}(0)^{^{\oplus b(V)}}
\oplus\cal O_{\!_\PP}(-1)^{^{\oplus c(V)}},
\quad a(V), b(V), c(V) \in \Bbb N \eqno(2.1)$$

Given $V\in\mm_C$ and $W\subset V|_{_C}$, a subbundle
with quotient isomorphic to ${\cal O}_{\!_\PP}(-1)$, put $V'={\cal
A}_C(V,W).$ Then $V'$ is
again in $\mm_C$, and we have
\[a(V') = a(V)+1\quad,\quad b(V') = b(V)\quad,\quad c(V') = c(V)-1\]
Therefore, the
operators $T_{{\cal O_{_{\!C}}(-1)}}$ and $T^*_{{\cal O_{_{\!C}}(-1)}}$
are well defined on the vector space
$\C[\mm_C]$, viewed as the subspace of functions on $\mm$ vanishing off
$\mm_C$.
We introduce the following modified operators
 $\Ee, \Ff$ on $\C[\mm_C]:$
$$(\Ee f)(V) = p^{-c(V)/2}\cdot(T^*_{_{{\cal O}_{_{\!C}}(-1)}}f)(V),\;\;
(\Ff f)(V) =  p^{-a(V)/2}\cdot(T_{_{{\cal O}_{_{\!C}}(-1)}}f)(V)
\eqno(2.2)$$
Define also the third operator $\Kk : \C[\mm_C] \to\C[\mm_C]$ by the
formula
\[ (\Kk f)(V) = p^{(a(V)-c(V))/2}\cdot f(V).\]

\ab Recall now, that to any generalized Cartan matrix one can associate
 a Kac-Moody Lie algebra $\g$. Furthermore, Drinfeld [Dr 1] and Jimbo [Ji] have
defined a $q$-deformation, ${\UU(\g)}$, of the universal enveloping
algebra of $\g$, the so-called quantized universal enveloping
algebra.

\ab The following crucial observation is at the origin of the
relationship
between quantum groups and Hecke operators. \hfill\break

{\sf Proposition 2.3.} {\it If the  curve $C$ is isomorphic to $\PP$
and has
self-intersection $(-2)$ in $S$, then
the operators $\Ee, \Ff$ and $\Kk$
satisfy
the commutation relations of the standard Chevalley generators of
the algebra $\bold U_{\bold q}(\goth s\goth l_2)$,
specialized at $\bold q=p^{1/2}$.}\hfill\break

{\sc Remark.} An analogue of the proposition still holds if the
projective line $C = \PP \subset S$ has an arbitrary  negative
self-intersection
$-m < 0$. Then, one should consider
vector bundles $V$ on $S$ such that $V|_{_C}$ is a direct sum
of several copies of the sheaves $\cal O_{\!_\PP}(0),\ldots,\cal
O_{\!_\PP}(m)\,.\,\square$

\ab Observe next that decomposition
of a vector
bundle $V$ on $\PP$ into a direct sum of line bundles, like (2.1), is not
 canonical in
general. What is only canonical is the filtration by subbundles
$0=\ldots\subseteq F_{j-1}V\subseteq
F_jV\subseteq F_{j+1}V\subseteq\ldots = V$
such that, for any $ j\in \Z$, the quotient
$F_{j}V/F_{j-1}V$ is the direct sum of several copies of
$\cal O_{\!_\PP}(-j)$.
The filtration $F_\bullet V$
measures the ``unstability" of the vector bundle $V$
and is called the Harder-Narasimhan filtration. Thus, for
any
point $x\in\PP$, the fiber $V_x$ of $V$ acquires a canonical finite
increasing filtration
$F_\bullet V_x$.

\ab We now turn to a more general case where
 $C\subset S$ is a reducible algebraic curve with
  smooth irreducible components $C_i, i\in \Delta$,
 isomorphic to $\PP$ each. Assume in addition that
 different components
intersect
each  other transversely, and
 each
$C_i$ has self-intersection $(-2)$. Thus, the
intersection matrix of the irreducible components of $C$ is the opposite
of a generalized Cartan matrix, $\|a(C)\|$.

\ab Call  a vector bundle $V$ on the curve $C$
{\it
admissible} if the following conditions hold :

\ab (i) \parbox[t]{106mm}{ For each irreducible component $C_i$, the
restriction
$V|_{_{C_i}}$ has a direct some decomposition like (2.1).}

\ab (ii) \parbox[t]{106mm}{ For any intersection point $x$ of two adjacent
components
$C_i$ and $C_j$, the two flags in the fiber at $x$ :
$\;(0=F_{-2}V_x\subseteq F_{-1}V_x\subseteq F_{0}V_x\subseteq
F_{1}V_x=V_x)$
and $(0=F'_{-2}V_x\subseteq F'_{-1}V_x\subseteq
F'_{0}V_x\subseteq F'_{1}V_x=V_x)\;$
induced by the canonical Harder-Narasimhan  filtrations on
$V|_{_{C_i}}$ and
$V|_{_{C_j}}$ are in generic relative position.}

\ab Fix an integer $d\geq 1$ and let $\mm_C$ be the moduli space of
$d$-dimensional
algebraic vector bundles on the surface $S$ that have admissible
restriction to $C$. For any $V\in \mm_C$,  $i\in \Delta$ and a
subbundle
$W\subset V|_{_{C_i}}$
such that $(V|_{_{C_i}})/W\simeq {\cal O_{C_i}(-1)}$,
the vector bundle ${\cal A}_{C_i}(V,W)$ is again
in $\mm_C$. Thus, the operators $T_{{\cal O}_{C_i}(-1)}$ and
$T^*_{{\cal O}_{C_i}(-1)}$
are well-defined on $\C[\mm_C]$ and,  for each $i\in\Delta$,
we introduce the operators
$\Ee_i, \Ff_i, \Kk_i ,$ defined as in (2.2), but with
respect to the corresponding component $C_i$. The following
result motivated by Nakajima [Na~1] gives way to an interpretation of
any quantized Kac-Moody algebra  $\bold U_{\bold q}(\g)$ in terms of
Hecke operators.
\hfill\break

{\sf Theorem 2.4.} {\it Let $\g$ be the Kac-Moody algebra whose Cartan
matrix is $\|a(C)\|$, the opposite of the intersection matrix.
Then, the operators}
$\{\Ee_i, \Ff_i, \Kk_i \in\;\mbox{End}\; \C[\mm_C],$
$i\in {\Delta}\}$ {\it satisfy the quantum Serre relations for the standard
Chevalley generators of the quantum algebra $\bold U_{\bold q}(\g)$,
specialized at $\bold q=p^{1/2}$.}\hfill\break

{\sc Remarks.}  (i) During the preparation of this paper we received
Nakajima's
preprint [Na 2] where a similar result is proved in the setup of
instantons on an
ALE-space (=desingularization of $\C^2/\Gamma$ for a finite subgroup
$\Gamma \subset SL_2(\Bbb C)$)
instead of algebraic vector bundles on a {\it compact} algebraic
surface.
Nakajima's proof is based on an
explicit
ADHM-type construction of instantons on an ALE-space in terms of a
linear algebra data. There is no ADHM-type construction for an arbitrary
surface, of course. Our approach is different. We
verify the Serre relations directly using a local geometry of vector
bundles on a small neighborhood of the curve $C\subset S$. Somewhat
miraculously,
the actual computation turns out to be quite similar to the
computation involving the Grassmannians of $k$-planes in $(\F_p)^d$,
carried out in [BLM].

\ab (ii)  There are many interesting example of configurations
of $\PP$'s on surfaces other than ALE-spaces. Take, for instance,
an Enriques surface
$S$ containing a configuration with intersection matrix
of affine type $E_8^{(1)}$, see [Ko]. The double cover $S'$ of $S$ is a
K3-surface
with two such configurations not intersecting each other.
Thus, we get an action on an appropriate space $\C[\mm_C]$ of two
commuting copies
of the quantum algebra $\bold U_{\bold q}(\g)$,
where $\g$ is the affine Lie algebra of type $E_8^{(1)}$. We would
like to thank R. Bezrukavnikov for pointing out this example.
$\quad\square$

\ab We now turn to Langlands reciprocity.
Let $S$ be an arbitrary smooth  projective
surface over $\Bbb F_p$, and $C\subset S$ an algebraic curve. Let
$\mm$ be the
set of $\Bbb F_p$-rational
algebraic
vector bundles on $S$ of rank $d$, and $\mm^\circ \subset \mm$ the
subset of all vector bundles that are trivial on $S \setminus C$.
Choose an integer $l$  prime to $p,$ and let
$\cal Loc_{_C}$ be the set of $l$-adic local systems
on $C$ of rank $d$.
Then we conjecture the following :

{\sf 2.5. Langlands conjecture  for surfaces.}
(i) {\it For any pair $(S,C)$, there is a natural algebra homomorphism}
\[\cal \C[\cal Loc_{_C}]\to \mbox{End}\, \C[\mm^\circ]\quad,\quad
f\mapsto T_f\]

(ii) {\it To any  $d$-dimensional local
system
$\phi$ on $S$, one can associate an ``automorphic" function
$F_\phi \in\C [\mm^\circ]$ such that we have
\[T_f(F_\phi) = f(\phi|_{_C})\cdot F_\phi\quad\mbox{for every}\quad
f\in \C[\cal Loc_{_C}].\]}

\ab Assume next that $S$ is a {\it complex} smooth projective surface with
a holomorphic symplectic
structure (e.g., a $K3$-surface) and $C\subset S$ a
connected algebraic curve.
Let $\mm$ denote the {\it moduli stack} of
algebraic (not necessarily stable)
vector bundles with trivial determinant on $S$, and let
$\mm^\circ$ denote the subspace of the vector bundles trivial on $S
\,\setminus\, C$. The symplectic
structure on $S$ induces a canonical symplectic structure on $\mm$.
The following conjecture is motivated,
in part, by the relationship, see Witten [Wi],
 between 2-dimensional field theory
and 3-dimensional Chern-Simons theory.

{\sf Conjecture 2.6.} $\mm^\circ$ {\it is a Lagrangian subvariety
in} $\mm$.}
\vskip 3mm
\pagebreak[3]
\ab {\bf 3. Quantum toroidal algebras.}

\ab  Given a complex semisimple Lie algebra $\g$, write $\g(z,u)$
for the
Lie algebra of all (Laurent) polynomial maps $\C^*\times\C^*\to\g$ with the
pointwise bracket. This `double-loop'
Lie algebra has a universal central extension by the infinite dimensional
vector space $\Omega^1/d\Omega^0$, where $\Omega^i$ stands for
polynomial
$i$-forms on $\C^*\times\C^*$. The extension is denoted
$\gt$ and  is called the toroidal Lie algebra associated to $\g$.
A description
of toroidal algebras in terms of generators and relations was first
given in [MRY].
We now sketch a construction of a
quantum deformation of the universal enveloping algebra
$U(\gt)$.

\ab Assume first that $\g=\goth s\goth l_2.$
The quantum toroidal algebra, $\UU(\goth s\goth l_{_{2,tor}})$, is a
$\C(q)$-algebra on
generators
\[E_0\,,\,E_1\,,\, E_2\;,\;F_0\,,\, F_1\,,\, F_2\;,\;K_0^{\pm 1}\,,\,
K_1^{\pm 1}\,,\,
 K_2^{\pm 1}, \]
that satisfy, in particular, the following relations
\[K_i\cdot K_i^{-1}=K_i^{-1}\cdot K_i=1\quad,\quad K_i\cdot K_j=K_j\cdot K_i,\]
\[K_i\cdot E_j\cdot K_i^{-1}=q^{a_{ij}}\cdot E_j\qquad,\qquad
K_i\cdot F_j\cdot K_i^{-1}=q^{-a_{ij}}\cdot F_j,\]
\[E_2\cdot E_0=q^{-2}\cdot E_0\cdot E_2\qquad,\qquad F_2\cdot
F_0=q^{-2}\cdot F_0\cdot
F_2,\]
involving the `double-extended' Cartan matrix:
$$
\|a_{i,j}\|\,=\,
\left(\matrix \;\;2&-2&\;\;2&\\
-2&\;\;2&-2\\
\;\;2&-2&\;\;2\endmatrix\right)\quad,\quad i,j = 0,1,2.$$
Further, if $|j-i|\leq 1$ we have :
\[[E_i,F_j]=\delta_{ij}\frac{K_i-K_i^{-1}}{q-q^{-1}}\]
and if $|j-i|=1$ the following Serre-type relations hold:
$$\sum_{m=0}^3(-1)^m\left[3\atop m\right]\cdot
E_i^{3-m}\cdot E_j\cdot E_i^m=0\quad,\quad
\sum_{m=0}^3(-1)^m\left[3\atop m\right]\cdot
F_i^{3-m}\cdot F_j\cdot F_i^m=0,$$
where $\left[i\atop j\right]$ denote
$q$-binomial coefficients. In addition, whenever $j-i=\pm 2$ we have
(put "$\pm$" respectively):
$$\sum_{m=0}^3(-q^2)^{\pm m}\left[3\atop m\right]\cdot
E_i^{3-m}\cdot F_j\cdot E_i^m=0\quad,\quad
\sum_{m=0}^3(-q^2)^{\pm m}\left[3\atop m\right]\cdot
F_i^{3-m}\cdot E_j\cdot F_i^m=0,$$
There are infinitely many other relations of higher order, however.
Therefore, for general semisimple Lie algebras,
 it will be more convenient to use an alternative
approach  similar to the ``loop-like" realization of quantum affine
algebras, introduced by Drinfeld [Dr 2].

\ab Thus, given a  semisimple Lie algebra $\g$,
let $\|a_{\al\be}\|_{(\al,\be=0,1,\ldots,r)}$ be the corresponding extended
(i.e., affine) Cartan matrix.
This matrix is symmetrizable,
that is, there are mutually prime integers $d_0,d_1,\ldots,d_r$ such that
the matrix with the entries $b_{\al\be}=d_\al\cdot a_{\al\be}$,
($\al,\be=0,1,\ldots,r$) is symmetric.
Thus, we have two
$(r+1)\times(r+1)$-matrices $\|a_{\al\be}\|$
and $\|b_{\al\be}\|$.

\ab Let $A$ be the free associative
$\C[q,q^{-1},$ $C,C^{-1}]$-algebra
on generators
$E_{\al,k}$, $F_{\al,k}$, $H_{\al,l}$, $K_{\al}^{\pm 1},$
($\alpha\in[0,r],\,k,l\in{\Bbb Z}, l\neq 0$).
Form the following generating functions :
\[K^{\pm}_\al(z) = K_\al^{\pm 1}\cdot \exp
\left((q - q^{-1}) \sum _{k \ge 1} H_\al(\pm k)\cdot z^{\mp
k}\right)\in A[[z^{\mp 1}]],\]
\[ E_\al(z) = \sum_{k = -\infty}^{\infty} E_{\al,k}\cdot z^{-k}
\enspace , \enspace
F_\al(z) = \sum_{k = -\infty}^{\infty} F_{\al,k}\cdot z^{-k}\enspace
\in A[[z,z^{-1}]],\]
\ab We define the quantum toroidal algebra, $\UU(\gt)$, to be the quotient
of the free algebra $A$ modulo a number of
relations. These are
most conveniently expressed in terms of the rational function
$\th_m(z)=\frac{q^m\cdot z -1}{z-q^m}$ and
the  generating functions $K^{\pm}_\al(z) , E_\al(z) , F_\al(z)$ above as
follows (cf. [Dr 2]) :
\[K^{\pm}_\al(z)\cdot K^{\pm}_\be(w)=K^{\pm}_\be(w)\cdot K^{\pm}_\al(z),\]
\[ \th_{b_{\al,\be}}(C^{-2} \cdot \zw) \cdot K_\al^{+}(z) \cdot
K_\be^{-}(w) = \th_{b_{\al,\be}}(C^{2}\cdot\zw)
\cdot K_\be^{-}(w) \cdot K_\al^{+}(z),\]
\[K_\al^{\pm}(z) \cdot E_\be(w) =
\th_{b_{\al,\be}}(C^{\pm 1} \cdot \zw) \cdot
E_\be(w) \cdot K_\al^{\pm}(z),\]
\[K_\al^{\pm}(z) \cdot F_\be(w) =
\th_{-b_{\al,\be}}(C^{\mp 1} \cdot \zw) \cdot
F_\be(w) \cdot K_\al^{\pm}(z),\]
\[ (q-q^{-1}) [E_\al(z)\,,\, F_\be(w)] = \de_{\al,\be}
\left( \de(C^{-2}\cdot \zw)\cdot K^{+}_\al(C\cdot w)
- \de(C^{2}\cdot \zw) \cdot K^{-}_\al(C\cdot z) \right),\]
\[ E_\al(z)\cdot E_\be(w) =
\th_{b_{\al,\be}}(\zw)\cdot E_\be(w) \cdot E_\al(z),\]
\[ F_\al(z)\cdot F_\be(w) =
\th_{-b_{\al,\be}}(\zw)\cdot F_\be(w) \cdot
F_\al(z),\]
and the symmetrizations with respect to $z_1,z_2,...,z_{1-a_{\alpha,\beta}}$
of each of the following two expressions vanish :
\[\sum_{i=0}^{1-a_{\al,\be}}(-1)^i
%\left[\begin{array}{c}1-a_{\al,\be}\\i\end{array}\right]
\left[{1-a_{\al,\be}}\atop i\right]
E_\al(z_1)\cdots E_\al(z_i) E_\be(w) E_\al(z_{i+1})\cdots
E_\al(z_{1-a_{\al,\be}})\]
\[\sum_{i=0}^{1-a_{\al,\be}}(-1)^i
%\left[\begin{array}{c}1-a_{\al,\be}\\i\end{array}\right]
\left[{1-a_{\al,\be}}\atop i\right]
F_\al(z_1)\cdots F_\al(z_i) F_\be(w) F_\al(z_{i+1})\cdots
F_\al(z_{1-a_{\al,\be}})\]

The specialization of the quantum toroidal algebra  at $q=1$
reduces
to the universal enveloping algebra $U(\gt)$. Furthermore, the
quantum toroidal algebra has the natural structure of a Hopf algebra.

\ab The double-loop Lie algebra $\g(z,w)$ has an obvious automorphism
switching the two variables $z$ and $w$. The following "quantized
analogue" of this observation, motivated by [Ch],
 clarifies the structure of $\UU(\gt)$.\hfill\break

{\sf Proposition 3.1.} {\it Let $U_1 \subset \UU(\gt)$ be the
subalgebra
generated by the elements $E_{\al,0}, F_{\al,0}, K_\al^\pm$,
$\al=0,1,...r$,
and let $U_2 \subset \UU(\gt)$ be the
subalgebra
generated by all the coefficients of the generating functions
$E_\al(z),F_\al(z), K_\al^\pm(z)$, $\al=1,2...r$ (note that $\al\neq
0$). Then we have

\ab (i) The algebras $U_1$ and $U_2$ are both isomorphic to
the quantum enveloping algebra with affine Cartan matrix
$\|a_{\al\be}\|_{\al,\be=0,...,r}$. Moreover, $U_1 \cap U_2=\UU(\g),$
is the quantum enveloping algebra with {\it finite} Cartan matrix
$\|a_{\al\be}\|_{\al,\be=1,...,r}$.

\ab (ii) There is an automorphism of the algebra $\UU(\gt)$ taking
$U_1$ to $U_2$.}

\vskip 3mm
\pagebreak[3]
\ab {\bf 4. Elliptic surface case.}

 \ab The considerations of section 2 become much more concrete for
elliptic surfaces. Thus, we fix in this section
a smooth projective surface $S$, a smooth curve $X$, and a
morphism $\pi : S \to X$ whose general fibre is an elliptic curve.
Let $o\in X$ be a point such that the fiber $C=\pi^{-1}(o)$
is {\it not} isomorphic to an elliptic curve. Such a fiber $C$
is known to be a connected rational curve. Call $C$ a {\it
non-multiple special fiber} if
the multiplicities of the irreducible components
of $C$ are relatively prime.
The local classification of non-multiple special fibers in elliptic surfaces,
due to Kodaira [Ko] and
N\'eron [N\'e], says that except a few cases,
they are all of the type
considered in section 2. More precisely, if we exclude
the special fibers of
type $A$, $C_1$, $C_2$ and $C_3$, as given in table [Ta, p.46],
the irreducible components of a non-multiple special fiber,
$C$, are isomorphic to $\PP$ each and
their intersection matrix is the opposite of
the Cartan matrix $\|a(C)\|$ of an affine root system.
Moreover, the correspondence $\|a(C)\| \leftrightarrow C$ sets up
a bijection between the affine root systems
of type $A^{(1)}, D^{(1)},E^{(1)}$ and the
possible types of non-multiple special fibers (different from
$A$, $C_1$, $C_2$ and $C_3$).
Let $\g_{_C}$ denote the affine Lie algebra associated to $\|a(C)\|$.
We thus obtain from theorem 2.4 the following result:\hfill\break

{\sf Theorem 4.1.} {\it The assignment
of the algebra generated by
the Hecke
operators}
$\{\Ee_i, \Ff_i, \Kk_i \in\;\mbox{End}\;\C[\mm_C]\,,\,
i\in \Delta\},$ {\it to a non-multiple special
fiber of an elliptic surface yields a bijective correspondence between
the possible types of a special fiber and the quantum affine algebras
of type $A^{(1)}, D^{(1)}$ and $E^{(1)}$.}$\;\square$

\ab It is likely that there is a similar interpretation of elliptic
algebras (see [GKV]) in terms of Hecke operators which correspond to a
{\it general}
non-singular fiber of an elliptic surface. When the fiber degenerates
to a special one, the elliptic algebra should degenerate to
a quantum affine algebra.

\ab We now extend theorem 2.4 to quantum toroidal algebras. To that
end, let $C^{\circ}\subset S$ be the curve obtained from a non-multiple
special fiber $C$ of the elliptic surface $S$ by removing one
geometric component
in such a way that the intersection matrix of the set
 $\Delta^\circ\subset\Delta$
of components of $C^\circ$ is of the corresponding (finite) type $A, D, E$. We
view
$C^{\circ}$ as a scheme (scheme structure is induced by that on $\pi^{-1}(o)$),
and let $n_i , i\in\Delta^\circ,$ denote
the multiplicities of the irreducible components of $C^{\circ}$.
Observe that the integers $n_i$
are the coordinates
of the maximal root of the root system $\Delta^\circ$ in the basis of
simple roots.

\ab In the setup of $\S 2$ we introduce two additional operators
$\Ee_0, \Ff_0$
corresponding to
the reducible curve $C^\circ$ as a whole by the formulas
$\Ee_0 = T_{_{{\cal O}_{_{C^\circ}}}}$,
$\Ff_0 = T^*_{_{{\cal O}_{_{C^\circ}}}}$. Recall $d=rk V,$ and put
\[(\Kk_0 f)(V) =
p^{^{(d-\sum_{i\in\Delta^\circ}n_i(a_i(V)-c_i(V)))/2}}
\cdot f(V).\]
\vskip 2mm
{\sf Theorem 4.2.} {\it Let $\gt$ be the toroidal algebra associated
with the affine Cartan
matrix $\|a(C)\|$, the intersection matrix.
Then, there is an algebra homomorphism} $\UU(\gt)| _{_{\bold q=p^{1/2}}}\to
\End\,\C[\mm_C]$.
{\it The image of this homomorphism is the algebra generated by the Hecke
operators}
$\{\Ee_i, \Ff_i, \Kk_i  \in\;\mbox{End}\;\C[\mm_C]\;,\,
i\in {\Delta}\cup\{0\}\}.$ \hfill\break

{\sf Remark 4.3.} The two subalgebras generated by the elements
$\{\Ee_i, \Ff_i, \Kk_i \,,\,
i\in {\Delta}\}$ and
$\{\Ee_i, \Ff_i, \Kk_i \,,\,
i\in {\Delta}^\circ\cup\{0\}\}$
 respectively correspond, via the theorem, to the
subalgebras $U_1$ and $U_2$ introduced in Proposition 3.1.$\quad\square$

\vskip 3mm
\pagebreak[3]
\ab {\bf 5. $\bold K$-theoretic construction of the toroidal algebra
for} $\bold {{\goth s\goth l}_n}$.

\ab Fix another positive integer $n$. In the special case $\g=\goth s\goth l_n$
the quantum toroidal algebra $\UU(\goth s\goth l_{_{n,tor}})$ has an
alternative definition. It is based on the well-known realization of the
affine Lie algera of type  $A_n^{(1)}$ as the subalgebra
of "$n$-periodic" elements in the affine Lie algebra
${\goth g}_\infty$ of type
$A_{\infty}^{(1)}$. Taking loops in both algebras
yields an imbedding of the toroidal Lie algebra
into ${\hat {\goth g}}_{\infty}$,
the central extension of the loop algebra on
${\goth g}_{\infty}$. This way, the
algebra $\UU(\goth s\goth l_{_{n,tor}})$ can be viewed as the
subalgebra of "$n$-periodic" elements of $\UU({\hat {\goth g}}_{\infty})$
the quantized enveloping algebra associated
(in the sense of [Dr 2]) with the Lie algebra
${\hat {\goth g}}_{\infty}$.

\ab It is useful to "enlarge" the algebra $\UU(\goth s\goth l_{_{n,tor}})$
slightly and introduce an extended  quantum toroidal algebra $\Uu$
with a central
extra-generator $D$ as follows.
Let $\|a_{\al\be}\|$ be the affine Cartan matrix of type $A^{(1)}_{n-1},$ and
$\|m_{\al\be}\|$ the $n\times n$-matrix given by
$m_{\al\be}=\delta_{\al,\be+1}-\delta_{\al+1,\be}$, where
$\al,\be\in\Z/n\Z$ and $\delta$ is
Kronecker's $\delta$.
Let $\Uu$ be quotient of the free associative
$\C[q,q^{-1},$ $C,C^{-1},D,D^{-1}]$-algebra
on the same set of generators as in $\S 3$
modulo the modified relations
\[K^{\pm}_\al(z)\cdot K^{\pm}_\be(w)=K^{\pm}_\be(w)\cdot K^{\pm}_\al(z),\]
\[ \th_{a_{\al,\be}}(C^{-2}D^{m_{\al\be}}\cdot \zw) \cdot K_\al^{+}(z) \cdot
K_\be^{-}(w) = \th_{a_{\al,\be}}(C^{2}D^{m_{\al\be}}\cdot\zw)
\cdot K_\be^{-}(w) \cdot K_\al^{+}(z),\]
\[K_\al^{\pm}(z) \cdot E_\be(w) =
\th_{a_{\al,\be}}(C^{\pm 1}D^{m_{\al\be}}\cdot \zw) \cdot
E_\be(w) \cdot K_\al^{\pm}(z),\]
\[K_\al^{\pm}(z) \cdot F_\be(w) =
\th_{-a_{\al,\be}}(C^{\mp 1}D^{m_{\al\be}}\cdot \zw) \cdot
F_\be(w) \cdot K_\al^{\pm}(z),\]
\[ (q-q^{-1}) [E_\al(z)\,,\, F_\be(w)] = \de_{\al,\be}
\left( \de(C^{-2}\cdot \zw)\cdot K^{+}_\al(C\cdot w)
- \de(C^{2}\cdot \zw) \cdot K^{-}_\al(C\cdot z) \right),\]
\[ E_\al(z)\cdot E_\be(w) =
\th_{a_{\al,\be}}(D^{m_{\al\be}}\zw)\cdot E_\be(w) \cdot E_\al(z),\]
\[ F_\al(z)\cdot F_\be(w) =
\th_{-a_{\al,\be}}(D^{m_{\al\be}}\zw)\cdot F_\be(w) \cdot
F_\al(z),\]
and same Serre-type relations as in $\S 3$.

\ab We produce below a geometric
construction
of the extended quantum toroidal algebra $\Uu$.
To that end, let $\K= \C((z))$ be the field
of Laurent formal power series. Given an integer $d\geq 1$,
write $\K^d$ for the $d$-dimensional coordinate vector
space over $\K$. By a lattice we mean a rank $d$ free
$\C[[z]]$-submodule
in $\K^d$.

\ab Following an idea of Beilinson-Lusztig-MacPherson [BLM]
in the finite case and its affine version introduced in [GV~1], let
$\bbb$ be the ``affine $n$-periodic partial flag variety" consisting
of the sequences
of lattices in $\K^d$ of the form
\[F= (\ldots\subseteq F_{i-1}\subseteq F_{i}\subseteq
F_{i+1}\subseteq\ldots)\enspace\mbox{such that}\enspace
F_{i+n}=z^{-1}\cdot F_i\enspace,\enspace\forall i\in \Z.\]
\ab Let $G(\K)=GL_d(\K)$ be the group of invertible $\K$-valued
$(d\times d)$-matrices. The natural $G(\K)$-action on the set
of lattices induces a $G(\K)$-action on $\bbb$ and the diagonal
$G(\K)$-action on $\bbb\times\bbb$. The set of $G(\K)$-orbits in
$\bbb\times\bbb$
is known [GV] to be in (1-1)-correspondence with the set of
 $\Z\times\Z$-matrices $\|a_{ij}\|$ subject to the following
conditions :

\ab (i) The entries of the matrix are non-negative integers such that,
for any $\hphantom{x}\qquad\qquad i,j\in\Z$, we have $a_{i+n,j+n}=a_{i,j}$;

\ab (ii) The matrix $\|a_{ij}\|$ has finitely many non-vanishing
diagonals, i.e., there exists $k>>0$ (depending on the matrix) such
that $a_{ij}=0$ whenever $|i-j|> k$;

\ab (iii) $\quad\quad\Sigma_{1\leq i\leq n}\Sigma_{j\in\Z} \;a_{ij} = d$

\ab The group $G(\K)$ is an
infinite-dimensional complex (algebraic) group. Similarly, the set $\bbb$
acquires
the natural structure, cf. [KL 1], of an infinite dimensional complex
variety, a direct limit of finite-dimensional projective varieties
of increasing dimension. Next, we define $T^*\bbb$, a kind of "cotangent
bundle" to $\bbb$ (which is not a vector bundle , cf. [Gi]), as follows.
Call an endomorphism $x\in {\goth g\goth l}_d(\K)$
{\it regular semisimple} (see [KL 3]) if it becomes regular semisimple after
an extension of scalars from $\K$ to ${\bar {\Bbb K}}$, an algebraic
closure of $\K$. Further, for any endomorphism $x\in
{\goth g\goth l}_d(\K)$ write its characteristic polynomial
\[\det(\lambda\cdot Id - x)= \lambda^d + \ldots \;\in\;
\K[\lambda]\;=\;\C((z))[\lambda]\]
Let ${\cal N_{rs}}$ be the set of all regular semisimple
endomorphisms of $\K^d$ such that its characteristic polynomial is of
the form $\det(\lambda\cdot Id - x)= \lambda^n\cdot P(\lambda,z)$
for some
$P \in \C[[z]][\lambda]$ (no negative powers of $z$ !). It is clear that
$0 \not\in {\cal N_{rs}}$ and that, for any
$x\in {\cal N_{rs}}$, there exists a lattice $L \subset \K^d$ such that $x(L)
\subset L$ and, moreover, $x^n(L)\subset z\cdot L$.
We define  the cotangent bundle
to $\bbb$ as follows
\[T^*\bbb = \{ (F,x) \in \bbb \times {\cal N_{rs}} \;|\;
x(F_i) \subset F_{i-1} \,, \forall i \in\Z \}\]
Dilation automorphisms $z\mapsto a\cdot z\,,\,a\in \C^*$ of the field
$\C((z))$
give rise to a natural ``rotation of the
loop"
$\C^*$-action on all the objects above.

\ab Following the pattern of [KL 2], [Gi], and [GV],  we introduce
a subvariety $Z\subset T^*(\bbb\times\bbb)$ to be the union of the
conormal bundles to the $G(\K)$-orbits in $\bbb\times\bbb$.
Clearly the variety $Z$ is $G(\K)$-stable with respect to the diagonal action.
Thus, $Z$ may be viewed as a
$(\C^*\ltimes G(\K))\times\C^*$-variety, where
the first copy of $\C^*$ accounts for the ``rotation of the
loop" action, and the last copy accounts for dilations along the
fibers of the cotangent bundle. Set $\G=(\C^*\ltimes G(\K))\times\C^*,$
and let $K^\G(Z)$ denote the Grothendieck group of $\G$-equivariant
coherent (in an appropriate sense) sheaves on $Z$.

\ab There is a convolution-type product on $K^\G(Z)$
defined  as follows (cf. [GV], [CG]).
Let $p_{ij} : T^*\bbb\times T^*\bbb\times T^*\bbb
\to$ $T^*\bbb\times T^*\bbb$ denote the projection along the factor
not named. Each  fiber of the map
\[ p_{13}: \enspace p_{12}^{-1}(Z) \cap p_{23}^{-1}(Z)\rightarrow
T^*\bbb\times T^*\bbb\]
were shown in [KL 3], to be a disjoint union of finite dimensional
projective varieties. Moreover, the image of this map equals $Z$.
Given $\G$-equivariant sheaves
${\cal F} , {\cal F'}$ on $Z$, define convolution
$K^\G(Z)\otimes K^\G(Z)\to K^\G(Z)$ by the formula
\[ [{\cal F}]\star[{\cal F'}]= (Rp_{13})_{\ast} \left( p_{12}^{\ast}
{\cal F} \bigotimes_{{}_{\Oc_{T^*\bbb\times T^*\bbb\times T^*\bbb}}}^{\L}
p_{23}^{\ast} {\cal F'} \right). \]
\vskip 3mm
\ab We have the following toroidal analogue
of [GV 1, theorem 7.9] :\hfill\break

{\sf Theorem 5.1.} {\it There is a surjective algebra homomorphism}
\[\Uu
\twoheadrightarrow K^{(\C^*\ltimes G(\K))\times\C^*}(Z).\]

All the homomorphisms above tend to become bijective as $d\to\infty$.
The central extra-generator $D$ in $\Uu$
goes under the morphism of theorem 5.1 to the generator of the
representation
ring $R(\C^*)\simeq \Bbb Z[D,D^{-1}]$, where the group
$\C^*$ acts via the ``rotation of the loop".
Theorem 4.2 and Theorem 5.1 may be viewed as the two components of an
analogue of the Satake isomorphism 1.1.

\vskip 1cm

\centerline{\bf References}

[BBD]$\enspace$ \parbox[t]{125mm}{
A. Beilinson, J. Bernstein, P. Deligne. Faisceaux Pervers. {\it Ast\'erisque}
 {\bf 100} (1981).}

[BLM]$\enspace$ \parbox[t]{125mm}{
 A. Beilinson, G. Lusztig, R. MacPherson. A geometric setting for quantum
groups. {\it Duke Math. J.} {\bf 61} (1990), 655-675.}

[Ch$\;$]$\enspace$ \parbox[t]{125mm}{
I. Cherednik. Double Affine Hecke
algebras, KZ-equations and Macdonald's operators.
{\it International Mathem. Research Notices, Duke Math. J.},
{\bf 9} (1992), 171-180.}

[CG$\;$] $\enspace$ \parbox[t]{125mm}{N. Chriss,
V. Ginzburg. Representation theory and Complex Geometry
(Geometric technique in Representation theory of Reductive groups).
{\it Progress in Mathem. Birkh\"auser} , 1995 (to appear).}

[Dr 1]$\enspace$ \parbox[t]{125mm}{
V. Drinfeld. Quantum Groups. {\it Proceedings of the ICM} , Berkeley 1986.}

[Dr 2]$\enspace$ \parbox[t]{125mm}{
V. Drinfeld. A new realization of Yangians and Quantum affine algebras.
{\it Soviet Math. Dokl.} {\bf 36} (1988), 212 - 216. }

[Dr 3]$\enspace$ \parbox[t]{125mm}{
V. Drinfeld. Two-dimensional $l$-adic representations of the
fundamental group of a curve over a finite field and automorphic forms
on $GL(2)$.
{\it Amer.J. Math.} {\bf 105} (1983), 85-114}

[Gi$\;$] $\enspace$ \parbox[t]{125mm}{
V. Ginzburg. Lagrangian construction of the universal enveloping
algebra
$U(sl_n)$. {\it C.R. Acad. Sci.} {\bf 312} (1991), 907-912.}

[GKV] $\enspace$ \parbox[t]{125mm}{
V. Ginzburg, M. Kapranov, and E. Vasserot. Elliptic algebras and
equivariant elliptic cohomology (to appear in {\it q-alg eprints}).}

[GV] $\enspace$ \parbox[t]{125mm}{
V. Ginzburg, E. Vasserot. Langlands Reciprocity for Affine Quantum
groups of type $A_n$, {\it Intern. Mathem. Research Notices. Duke
Math. J.}, {\bf 69}, 1993, p. 67-85.}

[IM$\;$]$\enspace$ \parbox[t]{125mm}{
N. Iwahori, H. Matzumoto. On some Bruhat decompositions and the structure
of the Hecke ring of a $p$-adic group. {\it Publ. Math\'em. I.H.E.S.}
{\bf 25} (1965), 5-48.}

[\;Ji$\;$\,]$\enspace$ \parbox[t]{125mm}{
M. Jimbo. A q-analogue of the enveloping algebra, Hecke algebra and
the Yang-Baxter equation. {\it Lett. in Mathem. Phys.}
{\bf  11} (1986), 247-252.}

[KL 1]$\enspace$ \parbox[t]{125mm}{
D. Kazhdan, G. Lusztig. Schubert varieties and Poincar\'e duality.
{\it Proc. Symp. Pure Math.} {\bf 36} (1980).}

[KL 2]$\enspace$ \parbox[t]{125mm}{
D. Kazhdan, G. Lusztig. Proof of the Deligne - Langlands conjecture for
affine Hecke algebras. {\it Invent. Math.} {\bf 87} (1987), 153-215.}

[KL 3]$\enspace$ \parbox[t]{125mm}{
D. Kazhdan, G. Lusztig. The fixed point variety in affine flag
manifold. {\it Israel J. Math.} {\bf 62} (1988), 129.}

[$\;$K$\;$] $\enspace$ \parbox[t]{125mm}{
K. Kodaira. On the structure of compact complex analytic surfaces.
{\it Amer. J. Math.}
{\bf 86} (1964), 751;
{\bf 88} (1966), 682; {\bf 90} (1969), 55.}

[$\;$Ko$\,$]$\enspace$ \parbox[t]{125mm}{
 S. Kondo. Enriques surfaces with finite automorphism
groups. {\it Japanese J. Math.} {\bf 12} (1986), 191-282.}

[$\;$L$\;$] $\enspace$ \parbox[t]{125mm}{
R. Langlands. Problems in the theory of automorphic forms.
{\it Lect. Notes in Math.} {\bf 170} (1970), 18-86.}

[$\;$La$\,$] $\enspace$ \parbox[t]{125mm}{
G. Laumon. Correspondance de Langlands G\'eom\'etrique pour les corps
de functions. {\it Duke Math. J.} {\bf 54} (1987), 309-359.}

[$\;$Lu$\,$] $\enspace$ \parbox[t]{125mm}{
G. Lusztig. Quivers, Perverse sheaves, and quantized enveloping
algebras.
{\it Journ. A.M.S.} {\bf 4} (1991), 365-421.}

[Ma$\;$] $\enspace$ \parbox[t]{125mm}{M. Maruyama. Elementary transformations
in the theory of algebraic vector bundles, in: ``{\it Algebraic
Geometry}, Lecture Notes in Math. {\bf 961}, Springer - Verlag, Berlin,
1983, p. 241-266. }

[MRY] $\enspace$ \parbox[t]{125mm}{
R.V. Moody, S.E. Rao, and T. Yokonuma. Toroidal Lie algebras and
Vertex representations. {\it Geometricae Dedicata} {\bf 35} (1990),
283-307.}

[Na 1]$\enspace$ \parbox[t]{125mm}{
H. Nakajima. Instantons on $ALE$-spaces, Quiver varieties, and
Kac-Moody algebras. {\it Duke Math. J.}  {\bf 76}, 2, (1994).}

[Na 2]$\enspace$ \parbox[t]{125mm}{
H. Nakajima. Gauge theory on Resolutions of Simple Singularities and
Simple Lie algebras. {\it Intern. Mathem. Research Notices.}
{\bf 2} (1994), 61-74.}

[N\'e$\;$] $\enspace$ \parbox[t]{125mm}{
A. N\'eron.  Mod\`eles minimaux des vari\'et\'es Ab\'eliennes sur les corps
locaux et globaux. {\it Publ. Math\'em. IHES.} {\bf 21} (1964), 5-126.}

[Ta$\;$]$\enspace$ \parbox[t]{125mm}{
J.T. Tate. Algorithm for determining the type of singular fiber in an elliptic
pencil. {\it Lecture Notes in Math.} {\bf 476} (1975), 33-52.}

[VW]$\enspace$ \parbox[t]{125mm}{
C. Vafa, E. Witten. A strong coupling test of $S$-duality.
Preprint , 1994.}

[Wi$\;$ ]$\enspace$ \parbox[t]{125mm}{
E. Witten. Topological Quantum Field theories.
{\it Commun. Math. Phys.} {\bf 117} (1988), 353.}

\vskip 10mm
\pagebreak[3]
V.G.: The University of Chicago, Mathem. Dept., Chicago, IL 60637, USA\\
$\hphantom{x} \hspace{10mm}$ginzburg@@math.uchicago.edu

M.K.: Northwestern University, Mathem. Dept., Evanston, IL 60208,
USA\\
$\hphantom{x} \hspace{10mm}$kapranov@@chow.math.nwu.edu

\nopagebreak
E.V.: Ecole Normale Sup\'erieure, 45 rue d'Ulm, 75005 Paris, FRANCE \\
$\hphantom{x} \hspace{10mm}$ vasserot@@dmi.ens.fr

\end{document}